\documentclass[11pt]{article}
\usepackage[a4paper,margin=27mm]{geometry}
\usepackage[UTF8,fontset=fandol,scheme=plain]{ctex}
\usepackage{amsmath,amssymb,amsthm,mathtools}
\usepackage{microtype,enumitem,comment}
\usepackage{tikz}
\usepackage{array,hhline}
\usepackage[font=small,labelfont=bf,format=hang]{caption}
\usepackage[colorlinks=true,linkcolor=blue,citecolor=blue,urlcolor=blue]{hyperref}
\numberwithin{equation}{section}
\newtheorem{theorem}{Theorem}

\newtheorem{lemma}[theorem]{Lemma}

\newtheorem{claim}{Claim}
\newcommand{\ex}{\operatorname{ex}}
\newcommand{\EX}{\operatorname{EX}}
\newcommand{\dist}{\operatorname{dist}}
\newcommand{\miss}{\operatorname{miss}}
\newcommand{\calL}{\mathcal L}
\newcommand{\calM}{\mathcal M}
\newcommand{\calP}{\mathcal P}
\newcommand{\calF}{\mathcal F}
\newcommand{\calS}{\mathcal S}
\newcommand{\eps}{\varepsilon}

\title{Approximating Tur\'an numbers by balanced product graphs}
\author{Chuandong Xu\thanks{School of Mathematics and Statistics, Xidian University, Xi'an, 710071, China. xuchuandong@xidian.edu.cn. Supported by the National Natural Science Foundation of China (No. 12571378).}}
\date{}

\begin{document}

\maketitle

\begin{abstract}
      Let $\calL$ be a finite family of graphs and put $p=\min_{L\in\calL}\chi(L)-1$. A graph is $\calL$-free if it contains no member of $\calL$ as a subgraph. The Tur\'an number $\ex(n,\calL)$ is the maximum number of edges in an $n$-vertex $\calL$-free graph. When $p\ge 2$, the Erd\H{o}s--Stone--Simonovits Theorem gives $\ex(n,\calL)=e(T_p(n))+o(n^2)$, where $T_p(n)$ is the balanced complete $p$-partite graph on $n$ vertices. The decomposition family $\calM$ of $\calL$ consists of the minimal bipartite graphs whose insertion into one part of $T_p(n)$ produces a graph containing a member of $\calL$. It is known that $\ex(n,\calM)=o(n^2)$.
      Suppose $p\ge 2$ and $\ex(n,\calM)=n^{\alpha+o(1)}$ for some $\alpha\in(1,2)$. We show that $\ex(n,\calL)$ differs by $o(\ex(n,\calM))$ from the maximum number of edges in an $\calL$-free graph containing $T_p(n)$ as a spanning subgraph.
\end{abstract}
\noindent\textbf{Keywords:} Tur\'{a}n number; extremal graph; product graph; forbidden scheme.

\section{Introduction}

      All graphs in this paper are simple. For a graph $G$, we write $v(G)$ for its number of vertices and $e(G)$ for its number of edges. A graph family $\calL$ is a finite set of graphs. We say that $G$ is \emph{$\calL$-free} if it contains no member of $\calL$ as a subgraph. The $\calL$-free graphs with the maximum number of edges for their order are called \emph{extremal graphs}. We write $\EX(n,\calL)$ for the set of such graphs on $n$ vertices and $\ex(n,\calL)$, the \emph{Tur\'an number} of $\calL$, for their number of edges. Let $T_p(n)$ be the \emph{Tur\'an graph}, the balanced complete $p$-partite graph
      on $n$ vertices with each part of size $\lfloor n/p\rfloor$ or
      $\lfloor n/p\rfloor+1$. We denote its number of edges by $t_p(n)$.
      
      In 1907, Mantel \cite{Mantel07} determined $\ex(n,K_3)$. Tur\'an \cite{Turan41} extended this result to complete graphs $K_{p+1}$ for all $p\ge 2$.

      \begin{theorem}[Tur\'an]\label{thm:Turan}
           For $p\ge 2$, $$\ex(n,K_{p+1})=t_p(n).$$
           Moreover, $T_p(n)$ is the only extremal graph.
      \end{theorem}

      Tur\'an proposed determining $\ex(n,L)$ for graphs other than complete graphs. This problem led to the study of extremal graph theory. The \emph{chromatic number} $\chi(L)$ of a graph $L$ is the smallest number of parts in a partition of $V(L)$ into independent sets. Put $p=p(\calL)=\min_{L\in\calL}\chi(L)-1$. The Erd\H{o}s--Stone--Simonovits Theorem \cite{ErdStone46,ErdSim66} estimates $\ex(n,\calL)$ in terms of the number of edges in $T_p(n)$.

      \begin{theorem}[Erd\H{o}s--Stone--Simonovits]\label{thm:ESS}
            Let $\calL$ be a finite family of graphs with $p(\calL)=p$. Then
            \[
                  \ex(n,\calL)=t_p(n)+o(n^2).
            \]
      \end{theorem}

      When the graph family $\calL$ contains a bipartite graph, we have $p(\calL)=1$, $t_p(n)=0$, and $\ex(n,\calL)=o(n^2)$. The problem of determining $\ex(n,\calL)$ when $p(\calL)=1$ is called the \emph{degenerate} Tur\'an problem; see F\"uredi and Simonovits \cite{FurSim13} for a survey. We consider the case $p(\calL)\ge 2$, in which every graph in $\calL$ is nonbipartite.

\subsection{Stability and product structure}
      For two graphs $G$ and $J$ on a common vertex set $V$, their \emph{edit distance} is
      \[
      \dist(G,J)=|E(G)\triangle E(J)|.
      \]
      In other words, $G$ can be obtained from $J$ by adding and deleting $\dist(G,J)$ edges.
      
      When $p(\calL)\ge 2$, the Asymptotic Structure Theorem of Erd\H{o}s and Simonovits \cite{Erdos68,Sim68Stability} states that for every $G\in \EX(n,\calL)$, there exists a $T_p(n)$ on $V(G)$ such that
      $
            \dist(G,T_p(n)) =o(n^2).
      $
      The following Stability Theorem due to Erd\H{o}s and Simonovits \cite{Erdos67,Sim68Stability} extends this conclusion to $\calL$-free graphs whose number of edges is close to $\ex(n,\calL)$.

      \begin{theorem}[The Stability Theorem] \label{thm:stability}
            Let $\calL$ be a graph family with $p(\calL)=p\ge 2$. If $G$ is an $n$-vertex $\calL$-free graph and $e(G) = \ex(n,\calL)- o(n^2)$,
            then there exists a $T_p(n)$ on $V(G)$ such that
            \[
                  \dist(G,T_p(n)) = o(n^2).
            \]
      \end{theorem}

      To describe the structure of extremal graphs more precisely, we need some definitions. In what follows, we assume that $n$ is sufficiently large. In particular, $n\ge l\cdot p(\calL)$, where $l=\max_{L\in\calL}v(L)$. We write $K_p(n_1,\ldots,n_p)$ for the complete $p$-partite graph whose parts have sizes $n_1,\ldots,n_p$. The \emph{disjoint union} of two graphs $J_1$ and $J_2$ is denoted by $J_1\sqcup J_2$. The \emph{product graph} $J_1\otimes\cdots\otimes J_p$ is obtained from the disjoint union of $J_1,\ldots,J_p$ by adding every edge between vertices in different $J_i$'s. Thus $K_p(n_1,\ldots,n_p)$, including $T_p(n)$ when the parts are balanced, is a product graph of edgeless graphs $\overline{K}_{n_i}$. 
      
      The \emph{decomposition family} $\calM$, short for $\calM(\calL)$, consists of the minimal bipartite graphs whose insertion into one part of $T_p(n)$ produces a graph containing a member of $\calL$. Equivalently, for a bipartite graph $M$ with at least one edge, we have
      \[
            M\in \calM \iff \exists\, L\in\calL,\ L\subseteq (M\sqcup\overline K_l)\otimes K_{p-1}(l,\ldots,l),\ \text{and $M$ is minimal}.
      \]
      For example, $\calM(K_p)=\{K_2\}$ when $p\ge 3$, and $\calM(K_3(2,2,2))=\{C_4\}$. 
      
      The next theorem is due to Simonovits \cite{Sim68Stability}. In the range $\ex(n,\calM)>cn$, the error term $O(\ex(n,\calM))$ can be sharp; see Simonovits \cite{Sim97Survey} for an example.

      \begin{theorem}[The Decomposition Theorem]\label{thm:decomposition}
            Let $\calL$ be a graph family with $p(\calL)=p\ge 2$ and decomposition family $\calM$. Then for every $G\in \EX(n,\calL)$, there exists a $K_p(n_1,\ldots,n_p)$ on $V(G)$ such that
            \[
                  \dist(G,K_p(n_1,\ldots,n_p)) =O(\ex(n,\calM)+n).
            \]
      \end{theorem}

      The Stability Theorem and the Decomposition Theorem describe the approximate structure of (near) extremal graphs. Much less is known about the exact structure of extremal graphs.  In many cases where the exact structure is known, the extremal graphs are product graphs. This is the case for edge-critical graphs (Simonovits \cite{Sim68Stability}) and for certain complete multipartite graphs (Tur\'an \cite{Turan41}, Erd\H{o}s and Simonovits \cite{ErdSim71}, Desai and Gupta \cite{DesGup25}). For further examples and references, see F\"uredi and Simonovits \cite[Section~2.1]{FurSim13}. 
      
      When $\calM$ contains no forest, the known extremal graphs were, for a long time, all product graphs. When $\calM$ contains a forest, however, an extremal graph for $\calL$ need not be a product graph; see Simonovits \cite{Sim82}, Griggs, Simonovits and Thomas \cite{GrigRubSim98}, and Peng, Song and Yuan \cite{PengSongYuan26} for some examples.

      Simonovits \cite{Sim97Survey} conjectured that when $\calM$ contains no forest, every extremal graph for $\calL$ is a product graph of $p(\calL)$ graphs. As he explained, this would reduce almost all nondegenerate extremal problems to degenerate ones. Recently, Xu \cite{Xu26A} found a counterexample to the conjecture. Thus, in general, we can expect better approximations of extremal graphs by product graphs, but not exact product structure. This motivates our study in this paper.
      
      For a $p$-partition $\calP=\{V_1,\ldots,V_p\}$ of $V(G)$, let $G_i=G[V_i]$. If $\calP$ maximizes the number of edges between different parts of $G$, we call it a \emph{maximum-cross-edge $p$-partition}. Our first result is a stability theorem: for any maximum-cross-edge $p$-partition of a near-extremal graph $G$, there is a product graph based on that partition and close to $G$ in edit distance.
 
      \begin{theorem}\label{thm:dist}
            Let $\calL$ be a graph family with $p(\calL)\ge 2$ and decomposition family $\calM$. Suppose $\ex(n,\calM)=n^{\alpha+o(1)}$, $\alpha\in(1,2)$. Let $G$ be an $n$-vertex $\calL$-free graph with $e(G)=\ex(n,\calL)-o(n^2)$, and let $\calP=\{V_1,\ldots,V_p\}$ be any maximum-cross-edge $p$-partition of $G$. Put $\beta=2/(3-\alpha)$. There is an $\calL$-free product $J'=G_1'\otimes\cdots\otimes G_p'$, with each $G_i'$ a spanning subgraph of $G[V_i]$, such that
            \begin{equation}\label{eq:dist}
                  \dist(G,J')=O\left(n^{\beta+o(1)}+\ex(n,\calL)-e(G)\right).
            \end{equation}
            In particular, if $\ex(n,\calL)-e(G)=o(\ex(n,\calM))$, then $\dist(G,J')=o(\ex(n,\calM))$. If $G\in\EX(n,\calL)$, then $\dist(G,J')=O(n^{\beta+o(1)})$.
      \end{theorem}

      The assumption $\ex(n,\calM)=n^{\alpha+o(1)}$ is reasonable because it is weaker than the conjectured asymptotic formula $\ex(n,\calM)= (c+o(1)) n^\alpha$ for some $c>0$ and rational $\alpha$ (Erd\H{o}s and Simonovits \cite[Conjecture~2]{ErdSim82Compact}; see also \cite[Conjecture~1.6]{FurSim13}).

\subsection{Forbidden schemes and the main theorem}

Let $\calL$ be a graph family with $p(\calL)=p\ge 2$, and let $L\in\calL$. Let $\Phi(L)$ be the set of all mappings from $V(L)$ to $[p]$. For $\phi\in\Phi(L)$ and $i\in[p]$, let $L_{\phi,i}$ be the subgraph of $L$ induced by $\{v\in V(L):\phi(v)=i\}$. The ordered list $L_{\phi,1},\ldots,L_{\phi,p}$ is the \emph{placement} of $L$ generated by $\phi$. There are $p^{|V(L)|}$ placements in total. 

Let $J=J_1\otimes\cdots\otimes J_p$ be a product graph in which $v(J_i)\ge l$. Since all cross-edges are present in $J$, every copy of $L$ in $J$ corresponds to a placement of $L$. We say that $J$ \emph{forbids} the placement $L_{\phi,1},\ldots,L_{\phi,p}$ of $L$ if $J_i$ is $L_{\phi,i}$-free for some $i\in[p]$. Thus
\[
      J_1\otimes \cdots \otimes J_p\ \text{is $\calL$-free} \iff \forall\, L\in \calL,\ \forall\, \phi\in \Phi(L),\ \exists\, i\in [p],\ J_i\ \text{is $L_{\phi,i}$-free}.
\]

A \emph{forbidden scheme} of $\calL$ is a mapping $\sigma:\{(L,\phi):L\in\calL,\ \phi\in\Phi(L)\}\to[p]$ such that $L_{\phi,\sigma(L,\phi)}$ has an edge for every $L\in\calL$ and every $\phi\in\Phi(L)$. Let $\calS(\calL)$ be the set of all forbidden schemes of $\calL$. To forbid a placement of $L$, it suffices for one factor $J_i$ to be $L_{\phi,i}$-free. The value $\sigma(L,\phi)$ specifies which factor must satisfy this condition. Let
\[
      \calL_i^\sigma:=\{L_{\phi,i}:L\in\calL,\ \phi\in\Phi(L),\ \sigma(L,\phi)=i\}.
\]
We have
\begin{equation}\label{eq:forbidden}
      J_1\otimes \cdots \otimes J_p\ \text{is $\calL$-free} \iff
            \exists\,\sigma\in\calS(\calL),\ \forall\, i\in [p],\ J_i\ \text{is $\calL_i^\sigma$-free}.
\end{equation}
The idea behind our forbidden schemes goes back to Simonovits \cite[Proposition~1]{Sim82}.

Each $M\in\calM$ corresponds to a placement of an $L\in\calL$, in which one $L_{\phi,i}$ is the disjoint union of $M$ and possibly some isolated vertices, while every $L_{\phi,j}$ with $j\ne i$ is edgeless. Therefore,
\begin{equation}\label{eqn:M_in_Lsigma}
      \text{$G$  is $\calL_i^\sigma$-free} \implies  \ \text{$G$ is $\calM$-free}.
\end{equation}

For a given forbidden family, we can simplify the description by grouping placements and retaining only the members of each $\calL_i^\sigma$ that are minimal under subgraph containment,
as in the definition of the decomposition family. Table~\ref{tab:octahedron} gives an example. We keep the present definitions because they are easier to use in the general setting. For extremal graphs that are product graphs, forbidden schemes reduce the description of their factors to separate extremal problems, as illustrated by the Octahedron Theorem of Erd\H{o}s and Simonovits \cite{ErdSim71}; see also \cite[Theorem 2.6]{FurSim13}.

\begin{table}[htbp]
  \centering
  \renewcommand{\arraystretch}{1.2}
  \begin{tabular}{w{c}{1.5cm}w{c}{0.5cm}w{c}{0.5cm}w{c}{1.5cm}w{c}{0.5cm}w{c}{1cm}w{c}{1cm}w{c}{1.5cm}w{c}{0.5cm}w{c}{1cm}w{c}{1cm}}
    \hhline{---~---~---}
    placements & & && $\sigma$ & $\calL^\sigma_1$ & $\calL^\sigma_2$ && $\sigma$ & $\calL^\sigma_1$ & $\calL^\sigma_2$ \\
    \hhline{---~---~---}
    type 1 & $C_4$ & $\emptyset$ && $\sigma_1$ & $\{C_4,P_3\}$ & $\{C_4\}$  &$\longrightarrow$& $\sigma_1$ & $\{P_3\}$ & $\{C_4\}$ \\
    type 2 & $\emptyset$ & $C_4$ && $\sigma_2$ & $\{C_4\}$  & $\{C_4,P_3\}$  && $\sigma_2$ & $\{C_4\}$ & $\{P_3\}$ \\
    type 3 & $P_3$ & $P_3$ &&  &   &   &&  & & \\
    \hhline{---~---~---}
  \end{tabular}
  \caption{The left block lists three types of modified placements for the octahedron $K_3(2,2,2)$. For example, $(C_4,\emptyset)$ means that $L_{\phi,1}$ contains a $C_4$, while the second entry imposes no condition on $L_{\phi,2}$. There are only $1\times 1\times 2=2$ modified forbidden schemes according to these three types. The families on the right retain only the minimal graphs in the corresponding families in the middle.}
  \label{tab:octahedron}
\end{table}

Now we state the main theorem of this paper. 

\begin{theorem}\label{thm:main}
      Let $\calL$ be a graph family with $p(\calL)=p\ge 2$ and decomposition family $\calM$, and let $\calS(\calL)$ be the set of forbidden schemes of $\calL$. Suppose $\ex(n,\calM)=n^{\alpha+o(1)}$, $\alpha\in(1,2)$. Let $J=J_1\otimes\cdots\otimes J_p$ have the maximum number of edges among the $n$-vertex $\calL$-free product graphs with $|J_i|\in\{\lfloor n/p\rfloor,\lfloor n/p\rfloor+1\}$ for each $i\in [p]$. Then
      \begin{equation}\label{eq:balance-product}
            e(J)=t_p(n)
            +\max_{\sigma\in\calS(\calL)}
                  \sum_{i=1}^{p}\ex\!\left(|J_i|,
                                    \calL_i^\sigma\right),
      \end{equation}
      and
      \begin{equation}\label{eq:main}
            \ex(n,\calL)= e(J) +o\!\left(\ex(n,\calM)\right).
      \end{equation}
\end{theorem}

It is easy to see that equation~\eqref{eq:balance-product} follows directly from \eqref{eq:forbidden}. To prove \eqref{eq:main}, we first establish Theorem~\ref{thm:dist}, which shows that every near-extremal $\calL$-free graph $G$ is close to a product graph $J'$. We then compare $J'$ with the balanced product graph $J$ using the same forbidden scheme. In Section~2 we give some lemmas that will be used. The proofs of Theorem \ref{thm:dist} and \ref{thm:main} are in Section~3.

\section{Preliminaries}

Let $G$ be a graph and let $\calP=\{V_1,\ldots,V_p\}$ be a $p$-partition of $V(G)$. For $v\in V(G)$, let $N_G(v)$ denote the neighborhood of $v$ in $G$. For $v\in V_i$, write
\[
 d_{\mathrm{int}}(v):=|N_G(v)\cap V_i|,\quad
 d_{\mathrm{cross}}(v):=|N_G(v)\setminus V_i|,\quad 
 d_{\mathrm{miss}}(v):=n-|V_i|-d_{\mathrm{cross}}(v).
\]
Intuitively, $d_{\mathrm{int}}(v)$ and $d_{\mathrm{cross}}(v)$ are the number of internal edges incident to $v$ in $G_i$ and the number of cross-edges incident to $v$, respectively. The total number of missing cross-edges of $G$ with respect to $\calP$ is
\begin{equation}\label{eq:missing-definition}
 \miss(G,\calP)=\sum_{i<j}\bigl(|V_i||V_j|-e_G(V_i,V_j)\bigr)
 =\frac12\sum_{v\in V(G)}d_{\mathrm{miss}}(v).
\end{equation}
Thus
\begin{equation}\label{eq:edge-count}
 e(G)=\sum_{i<j}|V_i||V_j|+\sum_i e(G_i)-\miss(G,\calP).
\end{equation}

The following lemma is a direct corollary of the Stability Theorem.

\begin{lemma} \label{lem:cutsize}
      Let $G$ be an $n$-vertex $\calL$-free graph with $e(G)=\ex(n,\calL)-o(n^2)$. Every maximum-cross-edge $p$-partition $\calP=\{V_1,\ldots,V_p\}$ satisfies
      \[
            \sum_i e(G_i)=o(n^2),\qquad \miss(G,\calP)=o(n^2),\qquad \text{and}\quad  
            |V_i|=n/p+o(n).
      \]
\end{lemma}

\begin{proof}
      Let $\calP'=\{V'_1,\ldots,V'_p\}$ be the $p$-partition of the copy of $T_p(n)$ on $V(G)$ given by the Stability Theorem. A maximum-cross-edge $p$-partition minimizes the number of edges inside the parts, so
      \[
            \sum_i e(G_i)\le \sum_i e(G[V'_i])=o(n^2).
      \]
      The Stability Theorem also gives $|e(G)-t_p(n)|\le\dist(G,T_p(n))=o(n^2)$. Since
      \[
            e(G)=\sum_{i<j}|V_i||V_j|-\miss(G,\calP)+\sum_i e(G_i) \quad \text{and}
                  \quad \sum_{i<j}|V_i||V_j|\le t_p(n),
      \]
      we obtain
      \[
            \miss(G,\calP)\le \sum_i e(G_i)+o(n^2)=o(n^2).
      \]

      On the other hand, 
      \(
            t_p(n)=e(G)+o(n^2)\leq \sum_{i<j}|V_i||V_j|+\sum_i e(G_i)+o(n^2).
      \)
      Thus
      \[
            t_p(n)-\sum_{i<j}|V_i||V_j|=o(n^2).
      \]
      Using $\sum_i|V_i|=n$ and $2t_p(n)=(1-1/p)n^2+O(1)$, we obtain
      \begin{align*}
            \sum_i (|V_i|-n/p)^2
                  &= \sum_i |V_i|^2- n^2/p\\
                  &= \Bigl(\sum_i |V_i|\Bigr)^2- 2\sum_{i<j}|V_i||V_j|-n^2/p\\
                  &=2t_p(n)-2\sum_{i<j}|V_i||V_j|+O(1)\\
                  &=o(n^2).
      \end{align*}
      Hence $|V_i|=n/p+o(n)$ for every $i$.
\end{proof}

The next two lemmas bound the number of vertices of large internal degree and the number of vertices of small degree in near-extremal $\calL$-free graphs. For related results on extremal graphs, see Erd\H{o}s and Simonovits \cite{ErdSim71}.

\begin{lemma}\label{lem:intdeg}
      Let $G$ be an $n$-vertex $\calL$-free graph with $e(G)=\ex(n,\calL)-o(n^2)$, and let $\calP=\{V_1,\ldots,V_p\}$ be any maximum-cross-edge $p$-partition of $G$. For every $\eps>0$, the number of vertices with $d_{\mathrm{int}}(v)\ge\eps n$ is $O(1)$.
\end{lemma}

\begin{proof}

Since $\calP$ is a maximum-cross-edge $p$-partition of $G$, we have $d_{V_j}(v)\ge d_{V_i}(v)$ for every $v\in V_i$ and every $j\ne i$. Let $A_i:=\{v\in V_i:d_{V_i}(v)\ge\eps n\}$ for $i\in [p]$. 

Suppose there exists $j\in[p]$ such that $A_j$ contains a set $W$ with $|W|\ge 2l/\eps^p$, where $l=\max_{L\in\calL}v(L)$. For each $\vec{x}=(x_1,\ldots,x_p)\in V_1\times\cdots\times V_p$,
let $d_W(\vec{x})$ be the number of vertices in $W$ adjacent to every vertex in the tuple. Each $w\in W$ has at least $\eps n$ neighbors in each part, so
\[
      \sum_{\vec{x}} d_W(\vec{x})\ge |W|(\eps n)^p.
\]
Count the pairs consisting of a tuple $\vec{x}$ and an $l$-subset of its neighbors $U\subseteq W$  in two ways. This gives
\[
      \sum_{U\in\binom{W}{l}}\prod_i \Bigl|V_i\cap\bigcap_{u\in U}N_G(u)\Bigl|= \sum_{\vec{x}}\binom{d_W(\vec{x})}{l}.
\]

Let $\bar{d}:=\dfrac{\sum_{\vec{x}}d_W(\vec{x})}{\prod_i|V_i|}\ge \dfrac{|W|\eps^p n^p}{n^p}\ge 2l.$ By convexity,
\[
      \sum_{\vec{x}}\binom{d_W(\vec{x})}{l} \ge \sum_{\vec{x}}\binom{\bar{d}}{l}\ge \sum_{\vec{x}}\frac{\left(\bar{d}/2\right)^l}{l!}\ge \frac{(|W|\cdot\eps^p\cdot n^p)^l}{2^l\cdot l!\cdot (n^p)^{l-1}}= \frac{\eps^{pl}}{2^l\cdot l!}|W|^l n^p.
\]

Hence some $l$-element subset $U$ of $W$ satisfies
\[
      \prod_i \Bigl|V_i\cap\bigcap_{u\in U}N_G(u)\Bigr|\ge cn^p
\]
for a constant $c>0$. In particular, each $B_i=V_i\cap\bigcap_{u\in U}N_G(u)$ has size at least $c'n$ for some constant $c'>0$. Each $B_i$ is disjoint from $U$. Choose an $l$-subset $U_i$  from each $B_i$ uniformly, independently across the parts. The expected number of missing cross-edges is at most 
\[
      l^2\cdot\miss(G,\calP)/(c'n)^2=o(1).
\]

Thus we can choose $U_1,\ldots,U_p$ so that all cross-edges between them are present. Every vertex of $U$ is adjacent to every vertex of each $U_i$, so $U,U_1,\ldots,U_p$ form a $K_{p+1}(l)$ in $G$, contradicting that $G$ is $\calL$-free. Therefore $|A_i|<2l/\eps^p$ for every $i\in[p]$.
\end{proof}

\begin{lemma}\label{lem:crossdeg}
      Fix $\delta>0$. Let $G$ be an $n$-vertex $\calL$-free graph with $e(G)=\ex(n,\calL)-o(n^2)$.
      For sufficiently large $n$, the number of vertices with degree less than $(1-1/p-\delta)n$ is at most
      \[
            2\cdot\frac{\ex(n,\calL)-e(G)}{\delta n}.
      \]
      In particular, each extremal graph has minimum degree at least $(1-1/p)n-o(n)$.
\end{lemma}

\begin{proof}
      Let $\calP=\{V_1,\ldots,V_p\}$ be a maximum-cross-edge $p$-partition of $G$. Lemma~\ref{lem:cutsize} gives $|V_i|=n/p+o(n)$ for each $i$. For every $v\in X:=\{v\in V(G):d_G(v)<(1-1/p-\delta)n\}$,
      \[
            d_{\mathrm{miss}}(v)\ge \left(n-\frac{n}{p}-o(n)\right)-\left(1-\frac{1}{p}-\delta\right)n=\bigl(\delta-o(1)\bigr)n.
      \]
      Fix a constant $c>0$ such that $lc<1/p$, and let $B=\{v\in V(G):d_{\mathrm{miss}}(v)>cn\}$. Since $\miss(G,\calP)=o(n^2)$, we have $|B|=o(n)$ and $|X|=o(n)$.

      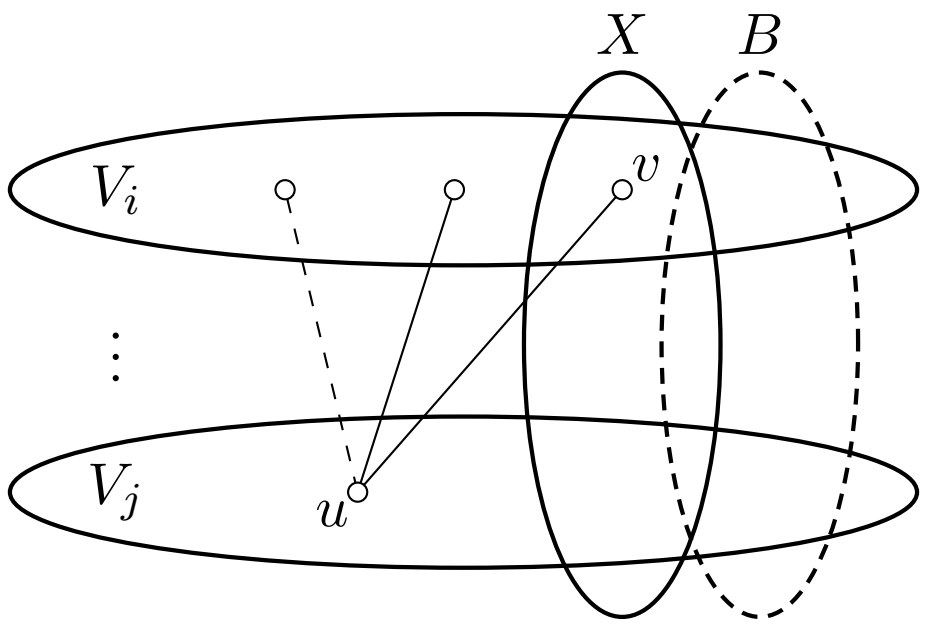
\begin{figure}[htbp]
            \centering
            \begin{tikzpicture}[x=1cm,y=1cm]
                  \foreach \y/\part in {3/i,1/j}{
                  \draw[thick] (0,\y) ellipse [x radius=3,y radius=0.5];
                  \node at (-2.30,\y) {$V_{\part}$};
                  }
                  \node at (-2.30,2) {$\vdots$};
                  \draw[thick] (1.05,1.975) ellipse [x radius=.65, y radius=1.80];
                  \draw[densely dashed,thick] (1.96,1.975) ellipse [x radius=.65,y radius=1.80];
                  \node at (1.05,4.03) {$X$};
                  \node at (1.96,4.03) {$B$};
                  \coordinate (vertex-v) at (1.05,3);
                  \coordinate (vertex-u) at (-0.70,1);
                  \coordinate (cross-nonneighbour) at (-1.18,3);
                  \coordinate (cross-neighbour) at (-0.06,3);
                  \draw (vertex-u) -- (vertex-v) (vertex-u) -- (cross-neighbour);
                  \draw[dashed] (vertex-u) -- (cross-nonneighbour);
                  \foreach \candidate in {cross-nonneighbour,cross-neighbour}{
                        \draw[fill=white] (\candidate) circle[radius=1.8pt];
                  }
                  \draw[fill=white] (vertex-v) circle[radius=1.8pt] node[above right,inner sep=1.5pt] {$v$};
                  \draw[fill=white] (vertex-u) circle[radius=1.8pt] node[below left,inner sep=1.5pt] {$u$};
            \end{tikzpicture}

            \caption{The partition of $G$ and the vertex subsets $X$ and $B$. The vertex $v$ is to be replaced in $V_i$. An already fixed neighbor $u$ of $v$ rules out at most $cn$ candidates.}
            \label{fig:intdeg}
      \end{figure}

      We construct an $\calL$-free graph $G'$ as follows:
      \begin{enumerate}[label=(\arabic*)]
            \item delete all edges incident to vertices in $X$;
            \item for each $i$, add all edges from $X\cap V_i$ to $V(G)\setminus(X\cup B\cup V_i)$;
            \item add all edges between $X\cap V_i$ and $X\cap V_j$ for $1\leq i<j\leq p$.
      \end{enumerate}
      Suppose $G'$ contains a copy $L^*$ of some $L\in\calL$. Keep the
      vertices of $L^*$ outside $X$ fixed, and replace those in $X$ one at
      a time by vertices in the same part outside $X\cup B$. After each replacement, keep the new vertex fixed and continue to denote the updated copy by $L^*$.
      Consider a vertex $v\in X\cap V_i$ of $L^*$ that is to be replaced.
      Every fixed neighbor of $v$ in $L^*$ lies outside $X\cup B\cup V_i$
      by construction. Thus each fixed neighbor of $v$
      has at most $cn$ nonneighbors in $V_i$ in $G$
      (see Figure~\ref{fig:intdeg}). Avoiding these nonneighbors and the
      vertices already used leaves at least $n/p-o(n)-lcn-l>0$ candidates. Thus we can finish the replacement of $v$. Edges between fixed vertices belong to $G$. At the end, we have a copy of $L$ in $G$, a contradiction.

      The number of edges deleted is at most $\sum_{x\in X}d_G(x)$, which is at most $|X|(1-1/p-\delta)n$. The construction adds at least $|X|[(1-1/p)n-o(n)]$ edges. For all sufficiently large $n$,
      \[
           \ex(n,\calL)-e(G)\ge e(G')-e(G) \ge \frac{1}{2}\delta |X|n.
      \]
      When $G\in\EX(n,\calL)$, we have $|X|=0$. Thus the minimum degree of $G$ is at least $(1-1/p)n-o(n)$.
\end{proof}

\begin{lemma}\label{lem:del_to_M-free}
      Let $\calL$ be a finite graph family with decomposition family $\calM$. Let $G$ be an $n$-vertex $\calL$-free graph with $e(G)=\ex(n,\calL)-o(n^2)$. Let $\calP=\{V_1,\ldots,V_p\}$ be any maximum-cross-edge $p$-partition of $G$, and let $G_i=G[V_i]$. After deleting $O(n+\ex(n,\calL)-e(G))$ internal edges, each $G_i$ becomes $\calM$-free. Moreover, both $\sum_i e(G_i)$ and $\miss(G,\calP)$ are
      \[
            O(\ex(n,\calM)+n+\ex(n,\calL)-e(G)).
      \]
\end{lemma}

\begin{proof}
      Fix two constants $\eps,\delta>0$ with $l(\eps+\delta)<1/p$. Let
      \[
            A:=\{v\in V(G):d_{\mathrm{int}}(v)\ge\eps n\},\quad X:=\{v\in V(G):d_G(v)<(1-1/p-\delta)n\}.
      \]
      Put $Z=A\cup X$. Lemmas~\ref{lem:intdeg} and \ref{lem:crossdeg} give $|Z|=O(1+(\ex(n,\calL)-e(G))/n)$. Note that $|V_i|=n/p+o(n)$ by Lemma \ref{lem:cutsize}. For every $v\in V_i\setminus Z$,
      \begin{equation}\label{eqn:missdeg}
            d_{\mathrm{miss}}(v)=n-|V_i|-d_G(v)+d_{\mathrm{int}}(v)\le(\eps+\delta+o(1))n.
      \end{equation}

      For $i\in[p]$, let $G_i^-$ be obtained by deleting from $G_i$ all edges incident to vertices in $Z$. The total number of deleted edges is $O(n+\ex(n,\calL)-e(G))$. Now we show that each $G_i^-$ is $\calM$-free. Suppose $G_i^-$ contains a copy of some $M\in\calM$. This copy avoids $Z$ because $M$ has no isolated vertices. 
      
      Use the placement of some $L\in\calL$ whose part $i$ is $M$ together with possible isolated vertices and whose other parts are edgeless. Keep the copy of $M$ and choose any additional vertices required in $V_i$ outside $Z$ and the copy. Then embed the vertices
      in the other parts one at a time, always outside $Z$. By
      \eqref{eqn:missdeg}, each previously embedded neighbor of the vertex to be embedded excludes at
      most $(\eps+\delta+o(1))n$ vertices. Avoiding these nonneighbors and the vertices already used leaves at least
      \[
            n/p-o(n)-l(\eps+\delta+o(1))n-l
      \]
      candidates. This lower bound is positive for sufficiently large $n$.
      Hence we can obtain a copy of $L$ in $G$, a contradiction. Thus
      \begin{align*}
            \sum_i e(G_i)&\le p\,\ex(n,\calM)+O(n+\ex(n,\calL)-e(G))\\
                  &= O(\ex(n,\calM)+n+\ex(n,\calL)-e(G)).
      \end{align*}
      
      The identity (\ref{eq:edge-count}) gives
      \begin{align*}
            \miss(G,\calP)&= \sum_{i<j} |V_i||V_j|-e(G) + \sum_i e(G_i)\\
                  &\le \ex(n,\calL)-e(G)+\sum_i e(G_i)\\
                  &= O(\ex(n,\calM)+n+\ex(n,\calL)-e(G)).
      \end{align*}
\end{proof}

When $G\in \EX(n,\calL)$, Lemma \ref{lem:del_to_M-free} implies that
\[
      \sum_i e(G_i) + \miss(G,\calP)=O(\ex(n,\calM)+n).
\]
This is the distance bound in the Decomposition Theorem. Moreover, Lemma \ref{lem:del_to_M-free} shows that the complete $p$-partite graph $K_p(n_1,\ldots,n_p)$ can be constructed from any maximum-cross-edge $p$-partition of $G$.

\begin{lemma}\label{lem:repair_to_J}
      Let $G$ be an $n$-vertex $\calL$-free graph with a partition $\calP=\{V_1,\ldots,V_p\}$ satisfying $|V_i|\ge cn$ for a fixed $c>0$ and all $i\in[p]$. Put $m:=\miss(G,\calP)$ and assume $m=o(n^2)$. Let $G_i=G[V_i]$. Suppose that after deleting $D_1(n)$ edges from $\bigcup_iG_i$, each remaining graph $G_i^-$ has maximum degree at most $D_2(n)$. Then after deleting $O(D_2(n)\sqrt m)$ further internal edges,
      the product $J'$ of the remaining graphs $G_i'$ is $\calL$-free and satisfies
      \[
            \dist(G,J')= O\bigl[D_1(n)+D_2(n)^2+\ex(n,\calL)-e(G)\bigr].
      \]
\end{lemma}

\begin{proof}
      If $m=0$ or $D_2(n)=0$, we can take $G'_i=G^-_i$. Then $\dist(G,J')=D_1(n)+m$. Assume from now on that $m\ge 1$ and $D_2(n)\ge 1$.

      Let $\calF$ be the family of all graphs with at most $l$ vertices and no isolated vertices, up to isomorphism. Put $u=2l^2D_2(n)\sqrt{m}$. In each $G_i^-$, whenever a graph $F\in\calF$ occurs and at most $u$ edges meet all its copies, delete such a set of edges. Continue until no such deletion is possible. Each deletion removes all copies of the chosen graph $F$ from that part, and later deletions cannot create new copies there. There are therefore at most $p|\calF|$ such deletions of sets of edges, removing at most $up|\calF|$ edges in total. Let $G_i'$ be the graph that remains in $V_i$.

      \begin{claim}
            The product graph $J'=G'_1\otimes\cdots\otimes G'_p$ is $\calL$-free.
      \end{claim}

      \begin{proof}
            Suppose not, i.e., the product graph contains a copy $L^*$ of a graph in $\calL$. Consider the subgraph $L^*[V_i]$ in each part. If it has an edge, let $L^-_i$ be the graph obtained by deleting all isolated vertices from $L^*[V_i]$. By the stopping condition, deleting all copies of $L^-_i$ in $G_i'$ takes more than $u$ edges. The edges of any maximal edge-disjoint packing meet all copies of $L^-_i$, so such a packing contains more than $u/e(L^-_i)$ copies. Choose a copy uniformly from this packing and then choose vertices uniformly from the remaining vertices of $V_i$ for the isolated vertices of $L^*[V_i]$. If $L^*[V_i]$ has no edges, choose all its vertices uniformly from $V_i$. When $L^*[V_i]$ has an edge, the probability of choosing a given vertex of $V_i$ is at most
            \[
                  \frac{D_2(n)}{u/e(L^-_i)}+\frac{l}{cn-l}\le\frac{l^2D_2(n)}{u}+\frac{l}{cn-l}.
            \]
            The first term uses the fact that a vertex lies in at most $D_2(n)$ copies in the edge-disjoint packing, since each copy containing it uses a distinct incident edge. The second term bounds the probability that the given vertex is chosen for one of the isolated vertices of $L^*[V_i]$. The same bound holds when $L^*[V_i]$ has no edges, where only the second term is needed. The choices in different parts are independent. The expected number of missing cross-edges between the chosen vertices is therefore at most
            \[
                  m\left(\frac{l^2D_2(n)}{u}+\frac{l}{cn-l}\right)^2
                     \le 2m\left(\frac{l^4D_2(n)^2}{u^2}+\frac{l^2}{(cn-l)^2}\right)
                     \le \frac12+\frac{2l^2m}{(cn-l)^2}<1
            \]
            for sufficiently large $n$. Thus we can choose the vertices so that all cross-edges between them are present. This gives a copy of $L^*$ in $G$, a contradiction.
      \end{proof}

      Let $D$ be the total number of internal edges deleted, including
      the initial deletions. Hence $D\le D_1(n)+up|\calF|$. 
      Since $J'$ is $\calL$-free and $e(J')=e(G)-D+m$, we have 
      \begin{align*}
            m&\le \ex(n,\calL)-e(G)+D\\
             &\le \ex(n,\calL)-e(G)+D_1(n)+2pl^2|\calF|D_2(n)\sqrt{m}\\
             &\leq \ex(n,\calL)-e(G)+ D_1(n) +  \frac{4p^2l^4|\calF|^2D_2(n)^2 +m}{2}.
      \end{align*}
      Rearranging gives
      $
            m/2\leq \ex(n,\calL)-e(G)+ D_1(n) + 4p^2l^4|\calF|^2D_2(n)^2.
      $
      Finally,
      \begin{align*}
            \dist(G,J')&\le D_1(n)+up|\calF|+m\\
                  &\leq D_1(n) + 2p^2l^4|\calF|^2D_2(n)^2 + \frac{m}{2} + m\\
                  &= O\bigl[D_1(n)+D_2(n)^2+\ex(n,\calL)-e(G)\bigr].
      \end{align*}
\end{proof}

For $2\le t\le n$, fix $G\in\EX(n,\calM)$ and choose a $t$-element subset of $V(G)$ uniformly. This induces at most $\ex(t,\calM)$ edges. Each edge belongs to the chosen set with probability $\binom t2/\binom n2$, so
      \[
            \ex(n,\calM)\frac{\binom t2}{\binom n2}\le\ex(t,\calM).
      \]
      Since $t^2/(2n^2)\le\binom t2\big/\binom n2$, we obtain
      \begin{equation}\label{eq:ex_le_ex}
            \frac{\ex(n,\calM)}{n^2}\le 2\frac{\ex(t,\calM)}{t^2}.
      \end{equation}
      In particular, $\ex(2t,\calM)\le8\ex(t,\calM)$.

\begin{lemma}\label{lemma:deg_sum}
For an $n$-vertex $\calM$-free
graph $G$, every $X\subseteq V(G)$ with $|X|\ge 2$ satisfies
\[
 \sum_{x\in X}d_G(x)\le 16\frac{n\ex(|X|,\calM)}{|X|}.
\]
\end{lemma}

\begin{proof}
      First suppose $|X|\le n/2$. Choose an $|X|$-subset $Y$ uniformly from $V(G)\setminus X$. Then
      \[
            \frac{|X|}{n-|X|}e_G(X,V(G)\setminus X)=\mathbb E[e_G(X,Y)]\le \ex(2|X|,\calM)\le 8\ex(|X|,\calM).
      \]
      Therefore,
      \begin{align*}
            \sum_{x\in X}d_G(x)&=2e(G[X])+e_G(X,V(G)\setminus X)\\      
                  &\le 2\ex(|X|,\calM)+8\frac{n-|X|}{|X|}\ex(|X|,\calM)\\
                  &\le 10\frac{n\ex(|X|,\calM)}{|X|}.
      \end{align*}
      If $|X|>n/2$, we have
      \[
            \sum_{x\in X}d_G(x)\le2\ex(n,\calM) \le2\ex(2|X|,\calM)
                  \le16\ex(|X|,\calM).
      \]
\end{proof}

\section{Proof of main results}

We first prove a distance bound in terms of a parameter $t$. A suitable choice of $t$ then gives Theorem~\ref{thm:dist}.

\begin{theorem} \label{thm:approx}
      Let $\calL$ be a finite graph family with $p=\min_{L\in\calL}\chi(L)-1\ge2$ and decomposition family $\calM$. There exists an integer $n_0\ge 2$ such that, for every fixed
      nonnegative function $r(n)=o(n^2)$, the following holds for all
      sufficiently large $n$ and every integer $n_0\le t\le n/4$.
      Let $G$ be an $n$-vertex $\calL$-free graph with
      $e(G)\ge \ex(n,\calL) - r(n)$.
      For any maximum-cross-edge $p$-partition $\calP=\{V_1,\ldots,V_p\}$ of $G$, there is an $\calL$-free product $J'=G_1'\otimes\cdots\otimes G_p'$, with each $G_i'$ a spanning subgraph of $G[V_i]$, such that
      \[
      \dist(G,J')= O\left[n+\ex(n,\calL)-e(G)+\frac{n\ex(t,\calM)}{t}
                        +\left(\frac{n\ex(t,\calM)}{t^2}\right)^2\right].
      \]
      The constant in the bound is independent of $G$, $\calP$, and $t$.
\end{theorem}

\begin{proof}
      By Lemma~\ref{lem:cutsize}, $|V_i|=n/p+o(n)$ for each $i$ and $\miss(G,\calP)=o(n^2)$. By Lemma~\ref{lem:del_to_M-free}, deleting $O(n+\ex(n,\calL)-e(G))$ edges from $\bigcup_iG_i$ leaves $\calM$-free graphs $G_i^-$. Fix $t\in[n_0,n/4]$ and put $T=40 n\ex(t,\calM)/t^2$.

      If $|V_i|\ge t$, let $X_i=\{x\in V_i:d_{G_i^-}(x)>T\}$. If $|X_i|\ge t$, choose
      $X\subseteq X_i$ with $|X|=t$. Lemma~\ref{lemma:deg_sum} gives
      \[
            tT<\sum_{x\in X}d_{G_i^-}(x)
            \le 16\frac{n\ex(t,\calM)}{t}\le tT,
      \]
      a contradiction. Thus $|X_i|<t$. Let $X_i^+$ be a $t$-subset of $V_i$ containing $X_i$, and delete all edges incident to $X_i$ from $G_i^-$. By Lemma~\ref{lemma:deg_sum} applied to $X_i^+$, we delete at most $16 n\ex(t,\calM)/t$ edges. The remaining graph has maximum degree at most $T$. If $|V_i|<t$, delete every edge of $G_i^-$; there are at most $\ex(t,\calM)$ such edges. The total number of deleted edges is
      \[
            D_1(n)=O\left(n+\ex(n,\calL)-e(G)+\frac{n\ex(t,\calM)}{t}\right).
      \]

      Applying Lemma~\ref{lem:repair_to_J} with $D_2(n)=T$ gives
      \begin{align*}
            \dist(G,J')&=O\bigl( D_1(n) + D_2(n)^2+\ex(n,\calL)-e(G)\bigr)\\
                  &=O\left[n+\ex(n,\calL)-e(G)+\frac{n\ex(t,\calM)}{t}
                        +\left(\frac{n\ex(t,\calM)}{t^2}\right)^2\right].
      \end{align*}
\end{proof}

\begin{proof}[Proof of Theorem \ref{thm:dist}]
      Take $t=\lceil n^{1/(3-\alpha)}\rceil$. Since $1<\alpha<2$,
      this lies in the range of Theorem~\ref{thm:approx} for sufficiently
      large $n$. Since $1+\dfrac{\alpha-1}{3-\alpha}=2+\dfrac{2\alpha-4}{3-\alpha}=\dfrac{2}{3-\alpha}=\beta$, we have
      \[
                  \frac{n\ex(t,\calM)}{t}=n^{\beta+o(1)},\qquad
                  \left(\frac{n\ex(t,\calM)}{t^2}\right)^2=n^{\beta+o(1)}.
      \]
      Since $1<\beta=\dfrac{2}{3-\alpha}<\alpha$, Theorem~\ref{thm:approx} gives 
      \[
            \dist(G,J')=O\!\left(n^{\beta+o(1)}+\ex(n,\calL)-e(G)\right).
      \]
      If $\ex(n,\calM)=\Theta(n^\alpha)$, the same choice of $t$ makes both terms involving $t$ of order $O(n^\beta)$. Furthermore, if $\ex(n,\calL)-e(G)=o(\ex(n,\calM))$, then
      \(
            \dist(G,J')=o(\ex(n,\calM)).
      \)
\end{proof}

To prove Theorem~\ref{thm:main}, we use two lemmas that quantify how balancing the part sizes affects the number of cross-edges and the extremal numbers of the factors.

\begin{lemma}\label{lem:sum-abs}
      Let $a_1\ge\cdots\ge a_p$ and $b_1\ge\cdots\ge b_p$ be nonnegative integers satisfying $\sum_i a_i=\sum_i b_i=n$ and $b_i\in\{\lfloor n/p\rfloor,\lfloor n/p\rfloor+1\}$ for each $i$. Then
      \[
            \Bigl(\sum_i|a_i-b_i|\Bigr)^2\leq2p\Bigl(\sum_{i<j}b_i b_j-\sum_{i<j}a_i a_j\Bigr).
      \]
\end{lemma}

\begin{proof}
      Let $r=n-p\lfloor n/p\rfloor$. If $r=0$, then
      $\sum_{i\le r}(a_i-b_i)=0$. If $r>0$ and
      $a_r\ge\lfloor n/p\rfloor+1$, then $a_i-b_i\ge0$ for every $i\le r$.
      If $r>0$ and $a_r\le\lfloor n/p\rfloor$, then
      $\sum_{i>r}(a_i-b_i)\le0$, so $\sum_{i\le r}(a_i-b_i)\ge0$.
      In all cases, we have
      \[
            \sum_i b_i(a_i-b_i)= \sum_{i\le r}(a_i-b_i) + \lfloor n/p\rfloor\sum_{i}(a_i-b_i) =\sum_{i\leq r} (a_i-b_i)\ge 0.
      \]
      Using $\sum_i a_i=\sum_i b_i=n$ again gives
      \[
            2\Bigl(\sum_{i<j}b_i b_j-\sum_{i<j}a_i a_j\Bigr)=\sum_i(a_i^2-b_i^2)=\sum_i(a_i-b_i)^2+2\sum_i b_i(a_i-b_i)\ge \sum_i(a_i-b_i)^2.
      \]

      By the Cauchy--Schwarz inequality,
      \[
            \Bigl(\sum_i|a_i-b_i|\Bigr)^2\le p\sum_i(a_i-b_i)^2\leq2p\Bigl(\sum_{i<j}b_i b_j-\sum_{i<j}a_i a_j\Bigr).
      \]
\end{proof}

\begin{lemma}\label{lem:ex-ex}
      Let $\sigma\in\calS(\calL)$ and $i\in[p]$, and put $\calF=\calL_i^\sigma$. If $a,b\in[n/(2p),n]$ are integers with $a,b\ge\max\{l,2\}$, then
      \begin{equation*}
            |\ex(b,\mathcal F)-\ex(a,\mathcal F)|\le 10p\,\frac{\ex(n,\calM)}{n}|b-a|.
      \end{equation*}
\end{lemma}

\begin{proof}
      Suppose $a\le b$. 
      Choose $a$ vertices uniformly from a graph in $\EX(b,\calF)$. Since the induced subgraph is $\calF$-free, its expected number of edges is at most $\ex(a,\calF)$. This gives
      \[
            \ex(b,\mathcal F) \dfrac{\binom{a}{2}}{\binom{b}{2}}\leq \ex(a,\mathcal F).
      \]
      Therefore,
      \begin{eqnarray*}
            \ex(b,\mathcal F)-\ex(a,\mathcal F)
                  &\le&\left[1-\frac{a(a-1)}{b(b-1)}\right]\ex(b,\mathcal F)\\
                  &=& \dfrac{(a+b-1)(b-a)}{b(b-1)}\ex(b,\mathcal F)\\
                  &\le& 10p\,\frac{\ex(n,\calM)}{n}|b-a|.
      \end{eqnarray*}
      The last inequality follows from $\ex(b,\mathcal F)\le \ex(n,\calM)$ (see (\ref{eqn:M_in_Lsigma})), $b\ge n/(2p)$, and $\dfrac{a+b-1}{b-1}\le 5$.
\end{proof}

\begin{proof}[Proof of Theorem \ref{thm:main}]
      Let $G\in\EX(n,\calL)$ and fix a maximum-cross-edge $p$-partition of $G$.
      Theorem~\ref{thm:dist} gives an $\calL$-free product $J'$ on this
      partition with $\dist(G,J')=O(n^{\beta+o(1)})$. By
      Lemma~\ref{lem:cutsize}, its factor sizes $a_i$ are $n/p+o(n)$. Let $b_i=|J_i|$.
      Order the parts so that $a_1\ge\cdots\ge a_p$ and
      $b_1\ge\cdots\ge b_p$. 

      Choose a forbidden scheme $\sigma$ for $J'$ by \eqref{eq:forbidden}.
      Realizing the same scheme on a balanced partition gives an
      $\calL$-free product graph with at most $e(J)$ edges. Thus
      \[
            e(J')\le \sum_{i<j}a_i a_j+\sum_i\ex(a_i,\calL_i^\sigma),
            \qquad
            e(J)\ge \sum_{i<j}b_i b_j+\sum_i\ex(b_i,\calL_i^\sigma).
      \]
      Lemmas~\ref{lem:sum-abs} and \ref{lem:ex-ex} give 
      \begin{align*}
            e(J') - e(J)&\le -\Bigl(\sum_{i<j}b_i b_j-\sum_{i<j}a_i a_j\Bigr)+\Bigl(\sum_i\ex(a_i,\calL_i^\sigma)-\sum_i\ex(b_i,\calL_i^\sigma)\Bigr)\\
                  &\leq - \Bigl(\sum_{i<j}b_i b_j-\sum_{i<j}a_i a_j\Bigr) + 10p\,\frac{\ex(n,\calM)}{n}\sum_i|a_i-b_i|\\
                  &\leq - \Bigl(\sum_{i<j}b_i b_j-\sum_{i<j}a_i a_j\Bigr) + 10p\,\frac{\ex(n,\calM)}{n}\sqrt{2p\Bigl(\sum_{i<j}b_i b_j-\sum_{i<j}a_i a_j\Bigr)}.
      \end{align*}

      Set $x=\dfrac{\ex(n,\calM)}{n}$ and $y=\sqrt{\sum_{i<j}b_i b_j-\sum_{i<j}a_i a_j}$. Then
      \[
            e(J') - e(J)\leq -y^2+10p\sqrt{2p}\,xy \le 50p^3x^2.
      \]
      Consequently,
      \[
            0\le\ex(n,\calL)-e(J)\le\dist(G,J')+50p^3x^2=O(n^{\beta+o(1)})+50p^3 \left[\frac{\ex(n,\calM)}{n}\right]^2.
      \]
      Note that $2\alpha-2<\alpha$ when $\alpha\in (1,2)$. Therefore, if $\ex(n,\calM)=n^{\alpha+o(1)}$, then
      \[
            e(G)-e(J)=o(\ex(n,\calM)).
      \]
\end{proof}

\section*{Statement of AI use}
The author proposed studying the error term in the
Erd\H{o}s--Stone--Simonovits Theorem.
ChatGPT obtained the initial version of Theorem~\ref{thm:dist},
which concerns approximation by $\calL$-free product graphs.
The author later suggested using placements of the forbidden graphs to describe the factors of these products. The definition of forbidden schemes emerged from discussions between the author and ChatGPT. The author rewrote the proofs
and assumes responsibility for all content.

\end{document}